\input ./style/arxiv-general.cfg
\documentclass[aap,MSNbibl,dvips]{arximspdf}
\makeatletter
   \@ifpackageloaded{graphicx}{}{\usepackage{graphicx}}
\makeatother

\doi{10.1214/15-AAP1121}
\volume{26}
\issue{3}
\pubyear{2016}
\firstpage{1942}
\lastpage{1946}
\docsubty{FLA}

\makeatletter
\newcommand{\eqref}[1]{(\ref{#1})}
\newproclaim{rem}{Remark}

\newcommand{\projected}[1]{|_{#1}}
\newcommand{\prob}{{\mathbb P}}
\newcommand{\tr}{\operatorname{tr}}
\newcommand{\R}{\mathbb{R}}
\makeatother

\begin{document}
\begin{frontmatter}

\title{A note on the expansion of the smallest eigenvalue distribution
of the LUE at the hard edge\thanksref{T1}}
\runtitle{Smallest eigenvalue expansion of LUE}

\begin{aug}
\author[A]{\fnms{Folkmar}~\snm{Bornemann}\corref{}\ead[label=e1]{bornemann@tum.de}}
\runauthor{F. Bornemann}
\affiliation{Technische Universit\"{a}t M\"{u}nchen}
\address[A]{Zentrum Mathematik---M3\\
Technische Universit\"at M\"unchen\\
80290 M\"{u}nchen\\
Germany\\
\printead{e1}}
\end{aug}
\thankstext{T1}{Supported by the DFG-Colla\-bo\-ra\-tive Research Center, TRR
109, ``Discretization in Geometry and Dynamics.''}

%
\received{\smonth{4} \syear{2015}}
%
\revised{\smonth{4} \syear{2015}}

%
\begin{abstract}
In a recent paper, Edelman, Guionnet and P\'{e}ch\'{e} conjectured a
particular $n^{-1}$
correction term of the smallest eigenvalue distribution of the Laguerre
unitary ensemble
(LUE) of order $n$ in the hard-edge scaling limit:
specifically, the derivative of the limit distribution, that is, the
density, shows up
in that correction term.
We give a short proof by modifying the hard-edge scaling to achieve an
optimal $O(n^{-2})$
rate of convergence of the smallest eigenvalue distribution. The
appearance of the
derivative follows then by a Taylor expansion of the less optimal,
standard hard-edge
scaling. We relate the $n^{-1}$ correction term further to the logarithmic
derivative of the Bessel kernel Fredholm determinant in the work of
Tracy and Widom.
\end{abstract}

%
\end{frontmatter}

\setcounter{footnote}{1}

\section{Introduction} We recall from random matrix theory that the
smallest eigenvalue distribution of the LUE of order $n$ with parameter
$a>-1$ is given by the Fredholm determinant \cite{ForrBook}, Chapter~9,
\[
\prob (\lambda_{\min} \geq s ) = \det \bigl(I- K_n^{a}
\projected{L^2(0,s)} \bigr),
\]
induced by the Laguerre projection kernel
\[
K_n^a(x,y) = \sum_{k=0}^{n-1}
\phi_k^a(x)\phi_k^a(y),\qquad\phi
_k^a(x) = \sqrt{\frac{k!}{\Gamma(k+a+1)}} e^{-x/2}
x^{a/2} L_k^{a}(x).
\]
Here, $L^a_k$ denotes the generalized Laguerre polynomial of degree $k$.
By Christoffel--Darboux and the relation $L^{a-1}_n(x) + L_{n-1}^{a}(x)
= L_n^{a}(x)$ (cf. \cite{MR0167642}, equation~(22.7.30)),
one gets the closed form
\[
K_n^a(x,y) = \frac{n!e^{-(x+y)/2}(xy)^{a/2}}{\Gamma(n+a)} \cdot\frac
{L_n^a(x)L_{n}^{a-1}(y)-L_{n}^{a-1}(x)L_n^a(y)}{x-y}.
\]
Using the Mehler--Heine type asymptotics (\cite{MR0372517}, Theorem~8.1.3) of the Laguerre polynomials, which holds uniformly
for bounded $z$ in the complex plane,\footnote{Hence,
for all $a\in\R$,
\[
z^{-a/2}J_a(2\sqrt{z}) = \sum_{k=0}^\infty
\frac{(-1)^k}{k! \Gamma
(a+k+1)}z^k
\]
constitutes a uniquely defined {\em entire} function of the complex
variable $z$.}
\[
n^{-a} L_n^{a}(z/n) = z^{-a/2}J_a(2
\sqrt{z}) + o(1) \qquad(n\to \infty)
\]
one immediately obtains, as first done by Forrester (\cite{MR1236195}, equation~(2.6); see also \cite{ForrBook}, Section 7.2.1), that in
the hard-edge scaling
\[
X = \frac{x}{4n},
\]
likewise for $Y$ and $y$, there is the limit, locally uniform for
positive $x$ and $y$,
\[
K_n^a(X,Y) \,dX = \bigl(K^a_\infty(x,y)
+ o(1) \bigr) \,dx \qquad (n\to\infty)
\]
with the Bessel kernel
\begin{eqnarray*}
K_\infty^a(x,y) &=& \frac{\sqrt{y} J_a(\sqrt{x})J_{a-1}(\sqrt
{y})-\sqrt{x} J_{a-1}(\sqrt{x})J_{a}(\sqrt{y})}{2(x-y)}
\\
&=& \frac{\sqrt{y} J_a(\sqrt{x}) J_a'(\sqrt{y}) -\sqrt{x}
J_a'(\sqrt{x}) J_a(\sqrt{y}) } {2(x-y)}.
\end{eqnarray*}
Lifted to the convergence of the induced Fredholm determinants (see the
next section for details),
one thus gets the hard-edge scaling limit of the LUE as a limit of
distributions, namely as $n\to\infty$
%
\begin{equation}
\label{eq:hard} F_n^a(s) = \prob \biggl(
\lambda_{\min} \geq\frac{s}{4n} \biggr) \to F_\infty^a(s)
= \det \bigl( I - K_\infty^a\projected {L^2(0,s)}
\bigr).
\end{equation}

Based on an identity of finite-dimensional Bessel function determinants
obtained from symbolical and numerical computer experiments, Edelman,
Guionnet and P\'{e}ch\'{e} \cite{1405.7590}, page~14, \emph{conjectured}
the following refinement of (\ref{eq:hard}):
%
\begin{equation}
\label{eq:conj} F_n^a(s) = F_\infty^a(s)
+ \frac{a}{2n} s f_\infty^a(s) + O\bigl(n^{-2}
\bigr),\qquad f_\infty^a(s) =\frac{d}{ds}F_\infty^a(s).
\end{equation}
In this note, we will give a short proof that this expansion, in fact,
holds true.

\begin{rem*} At FoCM'14, Gr\'{e}gory Schehr \cite{Schehr} presented
yet another proof (joint work with Anthony Perret) of
this expansion which
he had obtained as a byproduct of a new approach to the Painlev\'
{e}~III representation \cite{MR1266485} of the Bessel
kernel determinant.
\end{rem*}

\section{A short proof of the expansion} Lemma 4.1 of \cite
{1405.7590} easily implies a refinement of the Mehler--Heine type
asymptotics in the
form of an expansion (see also \cite{MR1509498}, page~29, for $a=0$ and
\cite{MR0029451}, page~156, for the
general case), namely,
%
\begin{eqnarray}
\label{eq:MehlerHeine}
&& (n+a)^{-a} L_n^{a}\bigl(z/(n+a)
\bigr)
\nonumber\\[-8pt]\\[-8pt]\nonumber
&&\qquad = z^{-a/2} J_a(2\sqrt{z}) -\frac{1}{2n}
z^{-(a-2)/2}J_{a-2}(2\sqrt {z})+ O\bigl(n^{-2}\bigr),
\end{eqnarray}
uniformly for bounded complex $z$. Subject to the following modified
hard-edge scaling, likewise for $Y$ and $y$,
\[
X = \frac{x}{4n} \biggl(1-\frac{a+c}{2n} \biggr),
\]
which transforms the Laguerre kernel according to
\[
K_n^a(X,Y) \,dX = \tilde K_n^a(x,y)
\,dx,
\]
this gives after some routine calculation\footnote{A {\em Mathematica}
notebook checking this result can be found at \arxivurl{arXiv:1504.00235}.}
the expansion
%
\begin{eqnarray}
\label{eq:kernelexp}
&& (xy)^{-a/2} \tilde K_n^a(x,y)
\nonumber\\[-8pt]\\[-8pt]\nonumber
&&\qquad = (x y)^{-a/2} K_\infty^a(x,y) -\frac
{c}{8n}
x^{-a/2} J_a(\sqrt{x}) y^{-a/2}J_a(
\sqrt{y}) + O\bigl(n^{-2}\bigr),
\end{eqnarray}
uniformly for bounded complex $x$ and $y$; in particular, uniformly for
$x,y \in[0,s]$. Because of the
pre-factor $(xy)^{-a/2},$ all the terms in \eqref{eq:kernelexp}
represent {\em entire} kernels; see the footnote in the last section.
Now, such a pre-multiplication of a kernel~$K$ by $(xy)^{-a/2}$ leaves
its Fredholm determinant invariant if we transform the measure defining
the underlying $L^2$-space accordingly:
\[
\det (I-K|_{L^2(0,s)} )= \det \bigl(I-(xy)^{-a/2}K|_{L^2((0,s);\nu_a)}
\bigr),\qquad\nu_a(dx) = x^{a} \,dx,
\]
which follows simply from conjugating $K$ with the \emph{unitary}
transformation
\[
U: L^2(0,s) \to L^2\bigl((0,s),\nu_a\bigr),
\qquad f(x) \mapsto x^{-a/2}f(x).
\]
Because of $a>-1$, the transformed measure $\nu_a$ on $[0,s]$ has
finite mass. Since the Fredholm determinant, when
defined with respect
to a measure $\nu$ of finite mass, is locally Lipschitz continuous
with respect to the \emph{uniform} convergence of the kernels (see
\cite{AGZ}, Lemma~3.4.5), the \emph{uniform} kernel expansion (\ref
{eq:kernelexp}) on $(0,s)^2$ immediately lifts to an expansion of the
induced Fredholm
determinants. For the particular choice $c=0$, which eliminates the
$O(n^{-1})$ term, we thus get
%
\begin{equation}
\label{eq:new} F_n^a \biggl( \biggl(1-\frac{a}{2n}
\biggr)s \biggr) = F_\infty^a(s) + O\bigl(n^{-2}
\bigr).
\end{equation}
Now, a simple Taylor expansion readily establishes \eqref{eq:conj}:
\[
F_n^a(s) = F_\infty^a \biggl(
\biggl(1-\frac{a}{2n} \biggr)^{-1}s \biggr) + O
\bigl(n^{-2}\bigr) = F_\infty^a(s) +
\frac{a}{2n} s f^a_\infty(s) + O\bigl(n^{-2}
\bigr).
\]
%
\begin{rem*}
The choice $c=-a$ of the additional scaling parameter which is
implicitly used in \eqref{eq:conj} is not the best
possible one; the optimally modified hard-edge scaling is given by
\eqref{eq:new}. Obtaining such a second-order convergence
rate by appropriately modifying the scaling was stimulated by the
corresponding work of Johnstone and Ma \cite{fast} for the \emph
{largest} eigenvalue distributions of the Gaussian unitary ensemble
(GUE) at the soft edge; cf. also the work
of Choup \cite{Choup} on the largest eigenvalue distributions of GUE
and LUE.
\end{rem*}

\section{Relation to the Tracy--Widom theory} Writing $\phi_a(x) =
J_a(\sqrt{x})$ for short and, as integral operators
acting on $L^2(0,s)$,
\begin{eqnarray*}
\hat K_n^a &=& K_\infty^a +
\frac{a}{8n}\phi_a\otimes\phi_a,
\\
I-\hat K_n^a &=& \bigl(I-K_\infty^a
\bigr) \biggl(I -\frac{a}{8n}\bigl(I-K_n^a
\bigr)^{-1} \phi_a \otimes\phi_a \biggr),
\end{eqnarray*}
we get from \eqref{eq:kernelexp} with $c=-a$, by the same reasoning as
in the last section,
\begin{eqnarray*}
F_n^a(s) &=& \det\bigl(I-\hat K_n^a
\bigr) + O\bigl(n^{-2}\bigr)
\\
&=& \det\bigl(I-K_\infty^a\bigr) \cdot\det \biggl(I -
\frac{a}{8n}\bigl(I-K_n^a\bigr)^{-1}
\phi_a \otimes\phi_a \biggr) + O\bigl(n^{-2}
\bigr)
\\
&=& F_\infty^a(s) \cdot \biggl(1- \frac{a}{8n}\bigl
\langle\bigl(I-K_n^a\bigr)^{-1}
\phi_a,\phi_a \bigr\rangle_{L^2(0,s)} \biggr) + O
\bigl(n^{-2}\bigr).
\end{eqnarray*}
Comparing the $n^{-1}$ term of the expansion with \eqref{eq:conj}
establishes the relation
%
\begin{equation}
\label{eq:u} -\frac{1}4 \bigl\langle\bigl(I-K_n^a
\bigr)^{-1} \phi_a,\phi_a \bigr
\rangle_{L^2(0,s)} = s \frac{f^a_\infty(s)}{F^a_\infty(s)} = s \frac{d}{ds} \log
F^a_\infty(s),
\end{equation}
which can already be found in the work of Tracy and Widom.
In fact, from the relation $\log\det= \tr\log,$ one gets that the
logarithmic derivative of
a Fredholm determinant $F(s)$ of a trace class operator $K$ acting on
$L^2(0,s)$ is generally given by \cite{MR1266485}, equation~(1.5),
\[
\frac{d}{ds}F(s) = \frac{d}{ds} \log\det\bigl(I-K
\projected{L^2(0,s)}\bigr) = -R(s,s),
\]
where $R(x,y)$ is the kernel of the operator $K(I-K)^{-1}$.
Specifically, in the case of the Bessel kernel, Tracy and Widom
calculated \cite{MR1266485}, equations~(2.5) and (2.21), that
\[
s R(s,s) = \tfrac{1}4 \bigl\langle\bigl(I-K_n^a
\bigr)^{-1} \phi_a,\phi_a \bigr
\rangle_{L^2(0,s)},
\]
which finally reproves \eqref{eq:u}.

\section*{Acknowledgement}
The author thanks Jamal Najim for suggesting to publish this note.




%

\printaddresses
\end{document}